\documentclass[12pt]{article}
\usepackage{amssymb,amsmath, amsthm}
\usepackage{color,xcolor}
\newtheorem{theorem}{Theorem}[section]
\newtheorem{lemma}[theorem]{Lemma}

\newtheorem{proposition}[theorem]{Proposition}

\newtheorem*{remark}{Remark}
\usepackage{mathtools}
\usepackage{array}
\usepackage{soul}
\date{}

\allowdisplaybreaks

\makeatletter
\renewcommand*\env@matrix[1][*\c@MaxMatrixCols c]{%
  \hskip -\arraycolsep
  \let\@ifnextchar\new@ifnextchar
  \array{#1}}
\makeatother

\begin{document}

\bibliographystyle{plain}

\title{Rigidity of proper holomorphic maps among generalized balls with Levi-degenerate boundaries}

\author{Sui-Chung Ng\footnote{School of Mathematical Sciences, East China Normal University, Shanghai, People's Republic of China. \textbf{Email:}~scng@math.ecnu.edu.cn}, Yuehuan Zhu\footnote{School of Mathematical Sciences, East China Normal University, Shanghai, People's Republic of China. \textbf{Email:}~52170601017@stu.ecnu.edu.cn}}

\maketitle

\begin{abstract}
   In this paper we studied a broader type  of generalized balls which are domains on the complex projective with possibly Levi-degenerate boundaries.	 We proved rigidity theorems for proper holomorphic mappings among them by exploring the structure of the moduli spaces of projective linear subspaces, which generalized some earlier results for the ordinary generalized balls with Levi-nondegenerate boundaries.
\end{abstract}

\section{Introduction}

We begin by extending the definition of generalized balls 
and let
\begin{align}\label{eq1}
\mathbb B^n_{r,s}= \biggr\{ [z_0,\dots,z_{n}]\in \mathbb{P}^{n}:\sum_{i=0}^{r}|z_i|^2 > \sum_{j=r+1}^{r+s}|z_i|^2 \biggr\},
\end{align}
for  nonnegative integers $ r,s $ with \(r+s \leq   n\). 
When   $ r+s =  n $,  we see that  the domain \(\mathbb B^n_{r,s}\) is precisely the ordinary generalized ball $ \mathbb B^n_r $ in $ \mathbb{P}^n $ defined by
$\mathbb B^n_r= \left\{[z_0,\dots,z_{n}]\in \mathbb{P}^{n}:\sum_{i=0}^{r}|z_i|^2 > \sum_{j=r+1}^{n}|z_i|^2 \right\}.$

The boundary $\partial\mathbb B^n_r$ of $\mathbb B_r^n$ is a homogeneous Levi-nondegenerate CR hypersurface in $\mathbb P^n$. The rigidity problems for the CR-maps defined on it and the closely related proper holomorphic maps on $\mathbb B^n_r$, have been extensively studied by many mathematicians, for instance, by Baouendi-Huang~\cite{baouendi2005super}, Baouendi-Ebenfelt-Huang~\cite{baouendi2011holomorphic}, Gao-Ng~\cite{gao2018rational}, Ng\cite{Ng1} and Seo~\cite{seo}. 

On the other hand, the boundary of $\mathbb B^n_{r,s}$ is Levi-degenerate, inhomogeneous and contains singularities when $r+s<n$. For the traditional methods in CR geometry using normal form theory, the degeneracy of the Levi-form causes a variety of difficulties, which in the case of $\partial\mathbb B^n_{r,s}$, cannot be overcome by taking higher derivatives of the local defining function. Nevertheless, there are works of Ebenfelt~\cite{ebenfelt}, Gaussier-Merker~\cite{gm}, Kossovskiy-Lamel-Xiao~\cite{klx}, etc. which dealt with Levi-degenerate CR hypersurfaces in various settings.  Specializing to those $\partial\mathbb B^n_{r,s}$ which are Levi-degenerate, Isaev-Kossovskiy~\cite{ik} also studied in detail their CR automorphisms.

In this paper, we investigate  the rigidity for proper holomorphic mappings among  $ \mathbb B^n_{r,s} $. Our first result is a generalization of a theorem of Baouendi-Huang~\cite{baouendi2005super} (Corollary 1.2 therein, for cases where $\ell\leq n/2$) regarding the rigidity of proper holomorphic maps between generalized balls. 
There is an interesting fact to note here. From Theorem~\ref{th3.2} below, we will see that for generalized balls with Levi-degenerate boundaries, it could happen that any local proper holomorphic maps between certain pair of such generalized balls must be \textit{totally degenerate}, in the sense that the differential of the map is everywhere degenerate. This is very different from the case for the generalized balls with Levi-nondegenerate boundaries.

 \begin{theorem}\label{th3.2}
 	Suppose \(1\leq r'\leq r < s \) and $U\subset\mathbb P^n$ is a connected open set such that $U\cap\partial\mathbb B^n_{r,s}\neq\emptyset$. If $f:U\rightarrow\mathbb P^m$ is a holomorphic map such that $f(U\cap \mathbb B^n_{r,s})\subset \mathbb B_{r'}^m$ and $f(U\cap \partial\mathbb B^n_{r,s})\subset \partial\mathbb B_{r'}^m$, then $r=r'$ and there exists $\Phi\in\mathrm{Aut}(\mathbb B^m_{r'})$ such that $\Phi\circ f([z_0,\ldots,z_n])=[z_0,\ldots,z_{r+s},0,\ldots,0]$.
 \end{theorem}

It is very natural to ask, after Theorem~\ref{th3.2}, what we can say about the local proper holomorphic maps from $\mathbb B^n_{r,s}$ to $\mathbb B_{r,s'}^{m}$. The first answer is that one can no longer conclude that such proper maps must be linear. This can be seen by considering the following simple example:
$$
[z_0,z_1,z_2,z_3,z_4]\mapsto [z_0^2,z_0z_1,z_0z_2,z_0z_3,z_1^2]
$$
which is easily seen to be a local proper holomorphic map from  $\mathbb{B}^{4}_{1,2}$ to $\mathbb{B}^{5}_{1,2}$.
Nevertheless, the next theorem, which may be regarded as a generalization of Theorem~\ref{th3.2}, provides a complete description of local proper holomorphic mappings from $\mathbb B^n_{r,s}$ to $\mathbb B_{r',s'}^{m}$ for $1\leq r'\leq r<s$. 

\begin{theorem}\label{th3.3}
Suppose \(1\leq r'\leq r < s \) and $U\subset\mathbb P^n$ is a connected open set such that $U\cap\partial\mathbb B^n_{r,s}\neq\emptyset$. If $f:U\rightarrow\mathbb P^m$ is a holomorphic map such that $f(U\cap \mathbb B^n_{r,s})\subset \mathbb B_{r',s'}^m$ and $f(U\cap \partial\mathbb B^n_{r,s})\subset \partial\mathbb B_{r',s'}^m$, then $r=r'$, $s\leq s'$ and there exist open set $U'\subset\mathbb P^n$ containing $[1,0,\ldots,0]$, $\Xi\in\mathrm{Aut}(\mathbb B^n_{r,s})$, $\Phi\in\mathrm{Aut}(\mathbb B^m_{r',s'})$ such that
$$
\Phi\circ f\circ\Xi(\left[1,z_1,\ldots,,z_n\right]) = [1,z_1,\ldots,z_{r+s},0,\dots,0, \psi_1, \ldots,\psi_{m-r'-s'}],
$$
on $U'$, where $\psi_1, \ldots,\psi_{m-r'-s'}$ are holomorphic functions on $U'$.
\end{theorem}

We now  give an   outline  of  the proofs of  Theorem~\ref{th3.2} and Theorem~\ref{th3.3}. Our basic strategy follows that of~\cite{Ng1}, but there are new technical obstacles to overcome. One difficulty arising here is that we can no longer assume that a proper holomorphic map is locally an embedding. (Indeed, as mentioned above, it turns out that the maps satisfying the hypotheses of Theorem~\ref{th3.2} are always totally degenerate if $r+s<n$.) Furthermore, when $r+s<n$, the boundary of $\mathbb B^n_{r,s}$ contains singularities that one needs to take care of. 

We begin by determining the sets of the maximal linear subspaces  contained in  $ \overline{\mathbb B}^{n}_{r,s}  $ and  $\mathbb B^n_{r,s}$. These linear subspaces turn out to be of dimension \(n-s\) and \(r\),   respectively. In addition, the former can be parameterized by  the points in the closure of the type-I bounded symmetric domain ${\Omega}_{r+1,s}$, while the latter can be parameterized by the  Cartesian product  of ${\Omega}_{r+1,s} $ and a Euclidean space.  This generalizes the result for the case $r+s=n$, which has been obtained in~\cite{Ng1}. Afterwards,  we  establish  that certain $r$-dimensional linear subspaces contained in the boundary of the source domain is mapped by the given proper holomorphic map to the maximal linear subspaces in the boundary of the target domain. From this and a maximum principle-type argument, we deduce that every $r$-dimensional linear subspace contained in $ {\mathbb{B}}^{n}_{r,s} $ is also mapped into some maximal linear subspace in  the target domain. This will then give the linearity of the map. 
For Theorem~\ref{th3.3}, after showing that the image of the boundary of the source domain will not lie completely inside in the singular part of $\partial\mathbb B^m_{r',s'}$, Theorem~\ref{th3.3} follows easily from Theorem~\ref{th3.2} once we compose the given proper map from $\mathbb B^n_{r,s}$ to $\mathbb B^m_{r',s'}$ with the canonical projection from $\mathbb{B}^{m}_{r',s'} $ onto $\mathbb B_{r'}^{r'+s'}$, which is also locally proper at smooth boundary points.

\section{On projective linear subspaces}\label{moduli}

At the beginning of this section, we recall some basic notions and facts  briefly.
Let $ M_{p,q}$  be the set of $ p\times q $ complex matrices. For a matrix $ X \in M_{p,q}$, the Hermitian transpose of $X$ is denoted by $ X^* $. If $ X $ is a Hermitian matrix, we denote by $ X >0  $ (resp. $ X \geq 0 $) the positive (resp. semi-positive )  definiteness of the matrix and write the $ p \times p $ identity matrix as $I_p$. Denote also by $ {\Omega}_{p,q} $ the classical bounded symmetric domain of type-I, which is defined as
\[ \Omega_{p,q}=\left\{Z \in M_{p,q}: I_p-Z Z^*>0\right\}. \]
Let  $G_{p,q}$ be  the Grassmannian of $p$-planes in  $\mathbb{C}^{p+q}$, which corresponds to the compact dual of $\Omega_{p,q}$. For a matrix \(X \in M_{p,p+q}\) of rank \(p\), we denote by $[X]\in G_{p,q}$  the $p$-plane spanned by the row vectors of  \(X\).  The matrix $X$ is uniquely determined by the point $ [X] \in   G_{p,q}$ up to  left multiplications by invertible matrices in $ M_{p,p}$. 
 Each matrix \(Z\in \Omega_{p,q}\) corresponds to the point \([I_p,Z] \in G_{p,q}\), which manifests the Borel embedding of \(\Omega_{p,q}\) into $G_{p,q}$. The generalized special unitary group $SU(p,q)$ is the subgroup of  the special linear group $SL(p+q,\mathbb{ C})$ consisting of matrices $X$ satisfying
\[ X\begin{bmatrix}
I_p &  0\\
0& -I_q
\end{bmatrix}X^*=\begin{bmatrix}
I_p &0  \\
0& -I_q
\end{bmatrix}. \]
In fact, the generalized ball $ \mathbb B^n_r $ is  an open orbit on $\mathbb{P}^{n}$ of $SU(r+1,n-r)$ under the standard action of $SL(n+1;\mathbb{C})$.

For a point $ [z]=[z_0,\cdots,z_n] \in \mathbb{P}^n  $, we split the homogeneous coordinates as $ [z]_{r,s}:=[z^1,z^2,z^3] $, where $ z^1=(z_0,\cdots,z_r) $, $ z^2=(z_{r+1},\cdots,z_{r+s}) $ and $ z^3=(z_{r+s+1},\cdots,z_n) $ (note that  if \(r+s=n\), \([z]\) is split as $ [z]_{r,s}:=[z^1,z^2] $). Let $\|z^1\|^2=\sum_{i=0}^{r} | z_i|^2$ and $ \|z^2\|^2=\sum_{j=r+1}^{r+s} | z_j|^2$. Then we have
\[
 \mathbb B^n_{r,s} = \left\{[z]_{r,s}\in \mathbb{P}^{n}: \|z^1\|^2 > \|z^2\|^2 \right\}.
 \]

\begin{lemma}\label{le1.3}
	 Let $ U'=\{[z]_{r,s}=[z^1,z^2,z^3]\in\mathbb{P}^n: (z^1,z^3)\neq 0\}$. If $ X\subset U' $ is a compact complex analytic subvariety of $ \mathbb{P}^n $, then $\dim(X)\leq n-s$. In particular, any compact complex analytic subvariety $X$ contained in $\overline{\mathbb B}_{r,s}^n$ is of dimension at most $n-s$.
\end{lemma}
\begin{proof}
	Consider the linear subspace $L^2:=\mathbb P^n\setminus U'=\{[z^1,z^2,z^3]\in\mathbb{P}^n: (z^1,z^3)=0\}$. Then $\dim(L^2)=s-1$. If $X\subset\mathbb P^n$ is a compact complex analytic subvariety such that $\dim(X)\geq n-s+1$, then $\dim(X)+\dim(L^2)\geq \dim(\mathbb P^n)$ and hence $X\cap L^2\neq\emptyset$. Thus, $X$ cannot be contained in $U'$. Finally, since $\overline{\mathbb B}_{r,s}^n\subset U'$, any compact complex analytic subvariety contained in $\overline{\mathbb B}_{r,s}^n$ is also of dimension at most $n-s$.
\end{proof}

\begin{proposition}\label{pro2.2}
	If  \(r \geq 1\), then the closure  $ \overline{\mathbb B}^{n}_{r,s}  $ contains a family of $(n-s)$-dimensional linear subspaces. They are compact  complex analytic subvarieties   of the maximal dimension in $ \overline{\mathbb B}^{n}_{r,s}  $. Furthermore, this set of $(n-s)$-dimensional  linear subspaces contained in  $ \overline{\mathbb B}^{n}_{r,s}  $ can be parameterized by the points in the closure $ \overline{\Omega}_{r+1,s} $.
\end{proposition}
\begin{proof} 
	For the case  \(r+s =   n\), the assertion has been proven by Ng~\cite{Ng1}.
	Thus it suffices to  prove the case where \(r+s <   n\).
We will give an explicit parameterization of a family of linear subspaces contained in  $ \overline{\mathbb B}^{n}_{r,s}  $.
Consider the $(n-s)$-dimensional linear subspace $L$  in $ \mathbb{P}^{n}  $ of the form
\begin{align}\label{eq3.1}
 L=\left\{[z]_{r,s}\in \mathbb{P}^{n}: z^2=(z^1,z^3)A \right\} \simeq \mathbb{P}^{n-s},
\end{align}
where   $ A \in M_{n-s+1,s}$ and  $ (z^1,z^3)=(z_0,\cdots,z_r,z_{r+s+1},\cdots,z_n) \in \mathbb{ C}^{n-s+1} $.  We have $L \subseteq \overline{\mathbb B}^{n}_{r,s} $ if and only if 
$$
   \left(z^{1}, z^{3}\right)\begin{bmatrix}
I_{r+1} &  0\\
0& 0
\end{bmatrix} \left(z^{1}, z^{3}\right)^{*} \geq \left(z^{1}, z^{3}\right) AA^* \left(z^{1}, z^{3}\right)^{*}   
$$
for all $ (z^1,z^3) \in  \mathbb{C}^{n-s+1}\setminus \{0\} $, which is in turn equivalent to
     $ \begin{bmatrix}
     I_{r+1} & 0 \\
     0& 0
     \end{bmatrix} - AA^*\geq 0 $. Here we partition the matrix $ A $  as the block matrix $  \begin{bmatrix}
     A_1   \\
      A_2
     \end{bmatrix} $, where $ A_1 \in M_{r+1,s}$ and  $ A_2 \in M_{n-r-s,s}$. Then 
     \[\begin{bmatrix}
     I_{r+1} & 0 \\
     0& 0
     \end{bmatrix} - AA^*=\begin{bmatrix}
     I_{r+1} & 0 \\
     0& 0
     \end{bmatrix}-\begin{bmatrix}
     A_1A_1^* & A_1A_2^* \\
    A_2A_1^* & A_2A_2^*
     \end{bmatrix} \geq 0
     \]
    is equivalent to $ I_{r+1} \geq A_1A_1^*  $ and  $ A_2=0 $. That is to say that the $(n-s)$-dimensional   linear subspace defined in Eq.(\ref{eq3.1}) lies inside $ \overline{\mathbb B}^{n}_{r,s}  $ if and only if $A=\begin{bmatrix}A_1\\0\end{bmatrix}$, in which $A_1 \in \overline{\Omega}_{r+1,s}$. 

    By Lemma~\ref{le1.3}, the $(n-s)$-dimensional linear  subspaces defined above  are of the maximal dimension among the compact complex analytic subvarieties contained in $ \overline{\mathbb B}^{n}_{r,s}  $. 
      \end{proof}

\begin{remark}
 When \(r+s <n\),  there is no  analogous one-to-one correspondence  between the  linear subspaces of the maximal dimension  in  $ {\mathbb{B}}^{n}_{r,s}  $ and the points of ${\Omega}_{r+1,s} $. In  fact, the dimension of the maximal
 linear subspaces contained in  $ {\mathbb{B}}^{n}_{r,s}  $ is less than \(n-s\).
  In addition, the  intersection of  $ {\mathbb{B}}^{n}_{r,s}  $ and
any linear subspace parameterized by the points of ${\Omega}_{r+1,s} $ in Proposition~\ref{pro2.2} is  a dense open subset of this given linear subspace, whose boundary points  lie on the boundary $ {\partial\mathbb{B}}^{n}_{r,s}   $.
\end{remark}

Next, we will study the maximal  linear subspaces  contained in $ {\partial\mathbb{B}}^{n}_{r,s}$ and their complex analytic properties.

\begin{lemma}\label{le2.3}
		 The compact complex analytic subvarieties of the maximal dimension contained in $ {\partial\mathbb{B}}^{n}_{r,s}   $ are  $({n-s})$-dimensional  linear subspaces for $r<s $ and  $({n-r-1})$-dimensional  linear subspaces for $ r\geq s $. Moreover,  if the image of any germ of a holomorphic map is contained in $ {\partial\mathbb{B}}^{n}_{r,s} $, then it must be contained in one of these linear subspaces.
\end{lemma}
\begin{proof}
	We only need to  give the proof for the case    \(r+s <   n\), because the   ordinary generalized ball situation has been proved by Ng~\cite{Ng1}.  Note that the boundary  $ {\partial\mathbb{B}}^{n}_{r,s}  \subset \mathbb{P}^{n}$ is defined as
	\[\left\lbrace  [z_0,\ldots,z_n]\in \mathbb{P}^{n}:|z_0|^2+\dots +|z_r|^2 = |z_{r+1}|^2+\dots +|z_{r+s}|^2 \right\rbrace.\]
	Suppose $r< s$. Let $ A $ = $  \begin{bmatrix}
	A_1   \\
	A_2
	\end{bmatrix} $, for $A_1 \in M_{r+1,s}$ and  $A_2 \in  M_{n-r-s,s}$. Then, by a similar argument as that in Proposition~\ref{pro2.2},
the $ (n-s)$-dimensional  linear subspace $ \left\{[z]_{r,s}\in \mathbb{P}^{n}: z^2=(z^1,z^3)A \right\} $  is contained in $ {\partial\mathbb{B}}^{n}_{r,s}$ if and only if $ A_1A_1^*=I_{r+1}   $ and $ A_2=0 $. These linear subspaces are also of the maximal dimension among the complex  linear subspaces of  $ {\partial\mathbb{B}}^{n}_{r,s} $ according to Proposition~\ref{pro2.2}. If $r \geq s $, by interchanging the coordinates $ z^1 $ with $ z^2 $ in $[z]_{r,s}$, we can also deduce that the  maximal linear subspaces contained in $ {\partial\mathbb{B}}^{n}_{r,s}   $ are    of dimension  ${n-r-1}$.
	
 The second part of the lemma will be proved in the same way as in~\cite{Ng1}.
	We will only prove the case $ r<s$ (and the case $r\geq s$ is similar). Let  $g: \bigtriangleup \rightarrow  {\partial\mathbb{B}}^{n}_{r,s} $ be a germ of a holomorphic curve, where $ \bigtriangleup$ is  the   complex unit disk. We may suppose that  $ g(\bigtriangleup) \subseteq {\partial\mathbb{B}}^{n}_{r,s} \bigcap U_0 $, where  $ U_0:=\left\lbrace[z] \in \mathbb{P}^{n}:z_0 \neq 0 \right\rbrace  $. We write $ g=(g_1,\cdots,g_{n}) $ in terms of the inhomogeneous coordinates in $U_0$, then $ g $ satisfies
	 \[ 1+\sum_{i=1}^{r}|g_i|^2 =\sum_{j=r+1}^{r+s}|g_j|^2. \]
	 Therefore, we deduce that  $ (g_{r+1},\cdots,g_{r+s})=(1,g_{1},\cdots,g_{r})H $ for some $ H \in M_{r+1,s}$ such that $ HH^*=I_{r+1}  $ (for a proof, see Lemma 4.1 in~\cite{ng2010holomorphic}). It now follows that such a germ of a holomorphic curve must be contained in one of those  $ (n-s) $-dimensional  linear subspaces  mentioned above.  Finally, the above argument can be applied similarly for any germ of  holomorphic map which sends the  unit ball in higher  dimension into ${\partial\mathbb{B}}^{n}_{r,s} $ and the proof is complete.
\end{proof}

We now look at the space of linear subspaces in $\mathbb B^n_{r,s}$. The space of $\ell$-dimensional linear subspaces in $\mathbb P^n$ is the Grassmannian $G_{\ell+1,n-\ell}$. As described before, for any point in $G_{\ell+1,n-\ell}$, we can find a matrix in $M_{\ell+1,n+1}$ of rank $\ell+1$ to represent it. If $A\in M_{\ell+1,n+1}$ is of rank $\ell+1$, then the corresponding $\ell$-dimensional linear subspace in $\mathbb P^n$ is just $L_A:=\{[zA]\in\mathbb P^n: [z]=[z_0,\ldots,z_\ell]\in\mathbb P^\ell\}$, in which $zA$ denotes the usual matrix multiplication between the row vector $z$ and the matrix $A$. The subset $E_\ell:=\{[A]\in G_{\ell+1,n-\ell}:[A]=[I_{\ell+1},A'] \,\textrm{\,for some\,}\, A'\in M_{\ell+1,n-\ell}\}$ is open and dense in $G_{\ell+1,n-\ell}$. Indeed, it is just a standard Euclidean cell in $G_{\ell+1,n-\ell}$ and it is easy to see that an $\ell$-dimensional linear subspace $L_A\subset\mathbb P^n$ defined above belongs to $E_\ell$ if and only if the projection $[z_0,\ldots,z_n]\mapsto [z_0,\ldots,z_\ell,0,\ldots,0]$ is everywhere defined on $L_A$. From the definition of $\mathbb B^n_{r,s}$, we see that the projection $[z_0,\ldots,z_n]\mapsto [z_0,\ldots,z_r,0,\ldots,0]$ is well defined on $\mathbb B^n_{r,s}$ and thus if an $r$-dimensional linear subspace $L_A$ is contained in $\mathbb B^n_{r,s}$, then $L_A\in E_r$. Thus, any $r$-dimensional linear subspace contained in $\mathbb B^n_{r,s}$ is of the form $\{[z,zA']\in\mathbb P^n:[z]\in\mathbb P^r\}$, where $A'\in M_{r+1,n-r}$. 

\begin{proposition} \label{pro2.5}
Suppose \(r+s <   n\). The set of \(r\)-dimensional linear subspaces contained in \(\mathbb B^n_{r,s}\) can be parameterized by  \(\Omega_{r+1,s} \times M_{r+1,n-r-s}\), where the latter is embedded in $G_{r+1,n-r}$ as an open set. Besides, these linear subspaces are  compact complex  analytic subvarieties   of the maximal dimension in \(\mathbb B^n_{r,s}\).
\end{proposition}
\begin{proof}
We just have remarked that any $r$-dimensional linear subspace contained in $\mathbb B^n_{r,s}$ is of the form
  \begin{equation}\label{eq4}
  	L_{A_1,A_2}= \left\{[z,zA_1,zA_2]_{r, s} \in \mathbb{P}^{n}: [z] \in \mathbb{P}^{r}\right\},
  \end{equation}
where $A_1\in M_{r+1,s}$, $A_2\in M_{r+1,n-r-s}$. Then, $ L_{A_1,A_2}$ is contained in \(\mathbb B^n_{r,s}\) if and only if \( \|z\|^2 > \|zA_1\|^2\) for all \([z] \in \mathbb{P}^{r}\), which is in turn equivalent to \(I-A_1A_1^*>0\).
Thus, the space of \(r\)-dimensional linear subspaces contained in \(\mathbb B^n_{r,s}\) can be parameterized by  \(\Omega_{r+1,s} \times M_{r+1,n-r-s}\).
  
  For the second part of  the proposition, the proof is similar to that of Lemma~\ref{le1.3}.
  \end{proof}

\begin{remark}
	If \(r+s<n\), then for each $l$ such that \(r<l\leq n-s\), the closure   $ \overline{\mathbb B}^{n}_{r,s}  $   contains a family of \(l\)-dimensional linear subspaces. However, 
all of	these  linear subspaces   will   intersect  the boundary ${\partial\mathbb{B}}^{n}_{r,s} $.  Moreover,  when \(r+s=n\),  the set of \(r\)-dimensional linear subspaces  in \(\mathbb B^n_{r,s}\)  \((\text{i.e.} \; \mathbb B^n_r)\; \)  is parameterized by the points of \(\Omega_{r+1,s}   \) (see Proposition 2.2 in~\cite{Ng1}).
	\end{remark}

The parametrization of the linear subspaces contained in $\mathbb B^n_{r,s}$ described in Proposition~\ref{pro2.5} can be extended to certain linear subspaces (of the same dimension) contained in $\partial\mathbb B^n_{r,s}$, as follows.

\begin{proposition} \label{pro2.6}
	If \(r<s\), there is a family of \(r\)-dimensional linear subspaces contained in \(\partial\mathbb B^n_{r,s}\), which can be parameterized by  $S(\Omega_{r+1,s}) \times M_{r+1,n-r-s}\subset G_{r+1,n-r}$ in the same way as in Eq.(\ref{eq4}) in Proposition~\ref{pro2.5}, where \(S({\Omega}_{r+1,s})\) is the Shilov boundary of \({\Omega}_{r+1,s}\). Moreover, the union of these linear subspaces is $\partial\mathbb B^n_{r,s}\setminus\{[z_0,\ldots,z_n]\in\mathbb P^n:z_0=\cdots=z_{r+s}=0\}$, which is just the smooth part of $\partial\mathbb B^n_{r,s}$.
\end{proposition}
\begin{proof}
By arguing in the same way as in Proposition~\ref{pro2.5}, we see that $L_{A_1,A_2}$ defined in Eq.(\ref{eq4}) is contained in $\partial\mathbb B^n_{r,s}$ if and only if $A_1\in S(\Omega_{r+1,s})$ and $A_2\in M_{r+1,n-r-s}$.

For the second part of the proposition, first of all it is easy to see that singular part of $\partial\mathbb B^n_{r,s}$ is $\{[z_0,\ldots,z_n]\in\mathbb P^n:z_0=\cdots=z_{r+s}=0\}$. Now let $p=[p_0,\ldots,p_n]$ be a point in the smooth part of $\partial\mathbb B^n_{r,s}$. Since $\sum^r_{j=0}|p_j|^2=\sum^{r+s}_{j=r+1}|p_j|^2$, we must have $(p_0,\ldots,p_r)\neq (0,\ldots, 0)$. Also, we can choose some $A_1\in M_{r+1,s}$ with $A_1A_1^*=I$ such that $(p_{r+1},\ldots,p_{r+s})=(p_0,\ldots,p_r)A_1$. Hence, $p\in L_{A_1,A_2}$ for this choice of $A_1\in S(\Omega_{r+1,s})$ and any $A_2\in M_{r+1,n-r-s}$.
\end{proof}

\section{Rigidity of proper holomorphic maps}

Now we are going to describe a double fibration pertaining to $\mathbb B^n_{r,s}$ and $\Omega_{r,s}$ which will be used to prove the rigidity of proper holomorphic maps between generalized balls. We begin with the  Grassmann bundles on the complex projective space $ \mathbb{P}^n $. For any point $ x \in \mathbb{P}^n $, denote by $ T_x \mathbb{P}^n$  its holomorphic tangent space at $ x $. The $ \ell $-Grassmann bundle $ G_\ell T \mathbb{P}^n $ $(1 \leq \ell \leq n)$ is given by the union of the Grassmannians $ G_\ell T_x \mathbb{P}^n $, i.e. $ G_\ell T \mathbb{P}^n:=\underset{x\in\mathbb{P}^n }{\bigcup }G_\ell T_x \mathbb{P}^n $, where the Grassmannian $ G_\ell T_x \mathbb{P}^n $ is the set of all  $\ell$-dimensional linear subspaces of $T_x \mathbb{P}^n$. Obviously, $ G_1T \mathbb{P}^n=\mathbb{P}(T \mathbb{P}^n) $ is just the projectivized tangent bundle of $ \mathbb{P}^n $.

For a fixed point \(x \in \mathbb{P}^n\), an element in $ G_\ell T_x \mathbb{P}^n $ is uniquely determined by a given $\ell$-dimensional linear subspace \(L \subseteq \mathbb{P}^n  \) with \(x \in \mathbb{P}^\ell \) (in general,  every pair  $(y, L)$ with \( y \in L\cong\mathbb{P}^\ell \subseteq \mathbb{P}^n \) will determine  an element in $G_\ell T_y \mathbb{P}^n$). Note that the Grassmannian \(G_{\ell+1,n-\ell}\) is the space which parameterizes all the \(\ell\)-dimensional linear subspaces of \( \mathbb{P}^n\). Then
it is clear that there is a holomorphic \(\mathbb{P}^\ell\)-foliation on $ G_\ell T \mathbb{P}^n $ and
$\pi^{\text{uni}} : G_\ell T \mathbb{P}^n \rightarrow G_{\ell+1,n-\ell} $ can be viewed as the universal family of all \(\ell\)-planes in \( \mathbb{P}^n\).  Therefore, we obtain a double fibration (cf.~\cite{ng2015holomorphic}) on $ G_\ell T \mathbb{P}^n $ as follows:
\begin{align} \label{eq3.5}
\mathbb{P}^n \stackrel{\pi}{\longleftarrow} G_\ell T \mathbb{P}^n \stackrel{\pi^{\text{uni}}}{\longrightarrow} G_{\ell+1,n-\ell}.
\end{align}

By considering the relation between the points in  $ \overline{\Omega}_{r+1,s} $ and the $(n-s)$-dimensional linear subspaces in $ \overline{\mathbb B}^{n}_{r,s}  $  described in Proposition~\ref{pro2.2}, we deduce  that  there is a closed subset \(E \cong \mathbb{P}^{n-s} \times \overline{\Omega}_{r+1,s}\) in  $G_{n-s}T \mathbb{P}^{n}$ such that the restriction of the universal family $\pi^{uni}:G_{n-s}T \mathbb{P}^{n} \rightarrow G_{n-s+1,s}$ to $E$ is just the canonical projection $\mathbb{P}^{n-s} \times \overline{\Omega}_{r+1,s}\rightarrow\overline{\Omega}_{r+1,s}$, and the fibers over \(\overline{\Omega}_{r+1,s}\) 
correspond to   those \((n-s)\)-dimensional linear subspaces  in $ \overline{\mathbb B}^{n}_{r,s}  $. As a matter of fact, if  \(r+s<n\),  $\pi^{uni}(E)$ is contained in some proper complex analytic  subvariety in $G_{n-s+1,s} $, in contrast to the situation where $r+s=n$ in which $\pi^{uni}(E)$ is just the closure of the open embedding $\Omega_{r+1,s}\subset G_{r+1,s}$, as demonstrated in~\cite{Ng1}.

Similarly, using Proposition~\ref{pro2.5}, we deduce that there is an open subset $W\cong \mathbb P^r\times \Omega_{r+1,s} \times M_{r+1,n-r-s}$ in $G_rT\mathbb{P}^n$ such that the restriction of the universal family $\pi^{uni}:G_rT \mathbb{P}^n\rightarrow G_{r+1,n-r}$ to $W$ is just the canonical projection $\mathbb P^r\times \Omega_{r+1,s} \times M_{r+1,n-r-s}\rightarrow \Omega_{r+1,s} \times M_{r+1,n-r-s}$. Moreover, $\pi^{uni}(W)\cong \Omega_{r+1,s} \times M_{r+1,n-r-s}$ is an open set in the standard Euclidean cell $\mathbb C^{(r+1)(n-r)}\cong M_{r+1,n-r}\cong M_{r+1,s}\times M_{r+1,n-r-s}$. On the other side of the double fibration described in Eq.(\ref{eq3.5}), we have $\pi(W)=\mathbb B^n_{r,s}$ since the fibers of $\mathbb P^r\times \Omega_{r+1,s} \times M_{r+1,n-r-s}\rightarrow \Omega_{r+1,s} \times M_{r+1,n-r-s}$ are just the $r$-dimensional linear subspaces in $\mathbb B^n_{r,s}$.

We are ready to prove our main theorems. 

\begin{proof}[Proof of Theorem~\ref{th3.2}]  
		In  case of the ordinary generalized ball $ \mathbb B^n_r $ (i.e. when $r+s=n$),  the assertion has been proven in~\cite{Ng1}. So 	we  only need to consider the case   \(r+s<n\).

Recalling the double fibration described above (Eq.(\ref{eq3.5})), we see that with respect to the universal family, the open subset $\pi^{uni}(G_{r}TU)\subset G_{r+1,n-r}$ corresponds to the set of those $r$-dimensional linear subspaces in $\mathbb P^n$ which intersect $U$. Let $S(\Omega_{r+1,s})$ be the Shilov boundary of $\Omega_{r+1,s}$. Since $U\cap\partial\mathbb B^n_{r,s}\neq\emptyset$ and by Proposition~\ref{pro2.6} the $r$-dimensional linear subspaces parametrized by $S(\Omega_{r+1,s})\times M_{r+1,n-r-s}\subset G_{r+1,n-r}$ are contained in $\partial\mathbb B^n_{r,s}$ and fill up the smooth part of $\partial\mathbb B^n_{r,s}$, it follows that $\pi^{uni}(G_{r}TU)$ intersects $S(\Omega_{r+1,s})\times M_{r+1,n-r-s}$. 

	 For any point $A\in \pi^{uni}(G_{r}TU)\cap (S({\Omega}_{r+1,s}) \times M_{r+1,n-r-s})$, we write $A=(A_1,A_2)$, where $A_1\in S({\Omega}_{r+1,s})\subset M_{r+1,s}$ and $A_2\in M_{r+1,n-r-s}$. Now we fix any such $A=(A_1,A_2)$ and choose a connected open set \(\mathcal V_1 \subseteq  M_{r+1,s}\cong\mathbb{C}^{(r+1)s}\) containing $A_1$ and a connected open set \(\mathcal V_2 \subseteq  M_{r+1,n-r-s}\cong\mathbb{C}^{(r+1)(n-r-s)}\) containing $A_2$ such that \(\mathcal V_1\times \mathcal V_2\subset \pi^{uni}(G_{r}TU)\). 
	 Take any point $z\in (\pi^{uni})^{-1}(A)\cap G_rTU$. We choose local inhomogeneous coordinates $(z_1,\ldots,z_r)$ on $(\pi^{uni})^{-1}(A)\cong\mathbb P^r$ such that $z=(0,\ldots,0)$ and choose also a small enough $r$-disk $\Delta^r\subset\mathbb P^r$ such that $\mathcal V_1\times\mathcal V_2\times\Delta^r$ is a local trivialization of the universal family bundle $G_rTU\rightarrow\pi^{uni}(G_rTU)$ at $(A,z)$.
	 
Now the composition $f^\flat:=f\circ\pi:\mathcal V_1\times\mathcal V_2\times\Delta^r:\rightarrow\mathbb P^m$ can be regarded as a holomorphic family of holomorphic maps from $\Delta^r$ to $\mathbb P^m$, where $\pi$ is the other canonical projection onto $\mathbb P^n$ in the double fibration Eq.(\ref{eq3.5}). In fact, it is evident from our construction that each member of this family is essentially the local restriction of $f$ to some $r$-dimensional linear subspace intersecting $U$. Recall that for any point $A_1\in \mathcal V_1\cap S(\Omega_{r+1,s})$, $A_2\in M_{r+1,n-r-s}$, the image $\pi((\pi^{uni})^{-1}((A_1,A_2)))$ is an \(r\)-dimensional linear subspace \(L_{A_1,A_2}=\{[z^1,z^2,z^3]_{r,s}\in \mathbb{P}^{n}: (z^2,z^3)=(z^1A_1,z^1A_2) \})\)  lying inside  $ {\partial\mathbb{B}}^{n}_{r,s}  $ (Proposition~\ref{pro2.6}). Now \(f\)  maps \(L_{A_1,A_2}\cap U\) into ${\partial\mathbb{B}}^{m}_{r}$ due to properness. Thus, \(f(L_{A_1,A_2}\cap U)\) is contained in an	 \(r'\)-dimensional linear subspace lying inside  $ {\partial\mathbb{B}}^{m}_{r'}  $  by  Lemma~\ref{le2.3}. That is, $f^\flat(A_1,A_2,\Delta^r)$ is contained in an $r'$-dimensional linear subspace for every $(A_1,A_2)\in  (\mathcal V_1\cap S(\Omega_{r+1,s}))\times\mathcal V_2$. 

By shrinking $\mathcal V_1,$ $\mathcal V_2$ and $\Delta^r$, we may assume that $f^\flat(\mathcal V_1\times\mathcal V_2\times\Delta^r)$ is contained in a Euclidean chart $\mathbb C^m\subset\mathbb P^m$ and thus we can write $f^\flat=(f^\flat_1,\ldots,f^\flat_m)$. We also write each $f^\flat_j$ as a power series  $f^\flat_j=\sum_{I}a_j^I(A_1,A_2)z^I$, in which $I$ is a multi-index,  $ z=(z_1,\ldots,z_r) \in \Delta^r $, and  $ a_j^I(A_1,A_2) $ are holomorphic functions on $\mathcal V_1\times\mathcal V_2$. Choose an arbitrary ordering $(I_0,I_1,\ldots)$ for $I$ with $I_0=(0,\ldots,0)$.
	Then, for a fixed $(A_1,A_2)$, the condition that $f^\flat(\{A_1\}\times\{A_2\}\times\Delta^r)$ is contained in an $r'$-dimensional linear subspace is equivalent to the condition that the rank of the following infinite matrix is at most $r'$:
	   \[ \begin{pmatrix}
	   	&a_1^{I_1}(A_1,A_2)  & a_2^{I_1}(A_1,A_2) & \ldots &a_{m}^{I_1}(A_1,A_2)  \\
	   	& \vdots & \vdots & \ddots &\vdots  \\
	   	&a_1^{I_k}(A_1,A_2)  & a_2^{I_k}(A_1,A_2) & \ldots &a_{m}^{I_k}(A_1,A_2)  \\
	   	& \vdots & \vdots & \ddots &\vdots  
	   \end{pmatrix}.
	   \]
	   
Since the rank condition can be rephrased as the vanishing of a set of minors of the above matrix, we see that for a fixed $A_2\in\mathcal V_2$, these minors (as holomorphic functions on $\mathcal V_1\times\{A_2\}$), vanish on $\mathcal V_1\cap S(\Omega_{r+1,s})$ due to the properness of $f$, as mentioned above. Then these minors (since they are holomorphic) actually vanish identically on $\mathcal V_1\times\{A_2\}$ (see~\cite{Ng1}, Lemma 2.9 therein for a proof). Since $A_2$ is also arbitrary, we conclude that for every $(A_1,A_2)\in\mathcal V_1\times\mathcal V_2$, $f^\flat(\{A_1\}\times\{A_2\}\times\Delta^r)$ is contained in an $r'$-dimensional linear subspace. Since $\pi^{uni}$ is open, we thus see that for a certain non-empty open subset $W\subset G_{r+1,n-r}$ and every $w\in W$, the image (under $f$) of the $r$-dimensional linear subspace given by $w$, i.e. $f(\pi((\pi^{uni})^{-1}(w))\cap U)$ is contained in some $r'$-dimensional linear subspace in $\mathbb P^m$. The containment in $r'$-dimensional linear subspace is clearly an analytic condition and thus we in fact have, for every $r$-dimensional linear subspace $L\subset\mathbb P^n$, $f(L\cap U)\subset L'$ for some $r'$-dimensional linear subspaces $L'\subset\mathbb P^m$. 

Now, since $r'\leq r$, by Lemma~\ref{le3.1}, either the entire image of $f$ is contained in an $r'$-dimensional linear subspace $L_0$ or $f$ extends to a linear rational map. Suppose $f$ is not linear and let $L_1\subset L_0$ be the linear span of the image of $f$. Thus, $d:=\dim(L_1)\leq r'$. Consider the intersection $L_1\cap\mathbb B^m_{r'}$. Since by hypotheses, $f(U\cap\mathbb B^n_{r,s})\subset\mathbb B^m_{r'}$ and $f(U\cap\partial\mathbb B^n_{r,s})\subset\partial\mathbb B^m_{r'}$, we see that under suitable coordinates on $L_1$, the intersection $L_1\cap\mathbb B^m_{r'}$ is also a generalized ball $\mathbb B_{a,b}^d$ for some $a,b\in\mathbb N$ and $f$ can be regarded as a holomorphic map $f^\star:U\subset\mathbb P^n\rightarrow\mathbb P^d$ such that $f^\star(U\cap\mathbb B^n_{r,s})\subset\mathbb B_{a,b}^d$ and $f^\star(U\cap\partial\mathbb B^n_{r,s})\subset\partial\mathbb B_{a,b}^d$. The maximal linear subspaces in $\partial\mathbb B_{a,b}^d$ are of dimension strictly less than $d$ and therefore if we repeat the entire argument on $f^\star$, then it follows that the image of $f^\star$ is contained in a proper linear subspace of $L_1$, contradicting the fact that $L_1$ is the linear span of the image of $f$. Thus, we conclude that $f$ extends to a linear rational map from $\hat f:\mathbb P^n\dasharrow\mathbb P^m$.

Finally, since $f(U\cap\partial\mathbb B^n_{r,s})\subset\partial\mathbb B^m_{r'}$, if we regard the linear map $\hat f$ as a linear map between two indefinite Hermitian vector spaces, then it maps an open subset of null vectors to null vectors. The set of null vectors is defined by a real quadratic equation and we thus see that $\hat f$ must map all null vectors to null vectors. Now, by Lemma~\ref{linear algebra}, it is a linear isometry up to a multiple $\lambda$. We must have  $\lambda>0$ since $f(U\cap\mathbb B^n_{r,s})\subset\mathbb B^m_{r'}$  by hypotheses (which translates to the condition that $\hat f$ maps some positive vectors to positive vectors). We then obtain that $r=r'$. By considering the image under $\hat f$ of the standard orthogonal basis, we see that for some  $\Phi\in\mathrm{Aut}(\mathbb B^m_{r'})$, we have $\Phi\circ f([z_0,\ldots,z_n])=[z_0,\ldots,z_{r+s},0,\ldots,0]$.
\end{proof}

\begin{proof}[Proof of  Theorem~\ref{th3.3}]
Let $\mathcal S=\{[z_0,\ldots,z_m]\in\mathbb P^m:z_0=\cdots=z_{r'+s'}=0\}$, which is the singular locus of $\partial\mathbb B^m_{r',s'}$. If $f(U\cap\partial\mathbb B^n_{r,s})\subset\mathcal S$, since $\partial\mathbb B^n_{r,s}$ is a real hypersurface in $\mathbb P^n$ and $\mathcal S$ is a complex linear subspace in $\mathbb P^m$, we must have $f(U)\subset\mathcal S\subset\partial\mathbb B^m_{r',s'}$, contradictory to the hypothesis that $f(U\cap\mathbb B^n_{r,s})\subset\mathbb B^m_{r',s'}$. Thus, $f$ at least maps a smooth point of $\partial\mathbb B^n_{r,s}$ to a smooth point of $\partial\mathbb B^m_{r',s'}$. Since the smooth part of the boundary is homogeneous (under the group of linear isometries preserving the associated Hermitian form), after replacing $f$ by its composition with some automorphisms of $\mathbb B^n_{r,s}$ and $\mathbb B^m_{r',s'}$, we may assume that $[1,0,\ldots,0]\in U$ and $f[1,0,\ldots,0]=[1,0,\ldots,0]$.

It is trivial that the canonical projection $\Pi:\mathbb B^m_{r',s'}\rightarrow\mathbb B_{r'}^{r'+s'}$, given by $[z_0,\ldots,z_m]\mapsto [z_0,\ldots,z_{r'+s'}]$, is locally proper in a neighborhood of any smooth point in $\partial\mathbb B^m_{r',s'}$ (in particular, near $[1,0,\ldots,0]$). Thus, the composition $\Pi\circ f:U\subset\mathbb P^n\rightarrow\mathbb P^{r'+s'}$ is a holomorphic map satisfying the hypotheses of Theorem~\ref{th3.2}. Therefore, $r=r'$, $s\leq s'$ and there exists $\phi\in\mathrm{Aut}(\mathbb B^{r'+s'}_{r'})$ such that $\phi\circ\Pi\circ f([z_0,\ldots,z_n])=[z_0,\ldots,z_{r+s},0,\ldots,0]$. Note that $\phi$ extends trivially to an automorhpism $\Phi$ of $\mathbb B^m_{r',s'}$ and we see that for some open subset $U'\subset U$, we have
$$
\Phi\circ f([1,z_1,\ldots,z_n])=[1,z_1,\ldots,z_{r+s},0,\dots,0, \psi_1, \ldots,\psi_{m-r'-s'}]
$$
for some holomorphic functions $\psi_1, \ldots,\psi_{m-r'-s'}$ on $U'$.
\end{proof}

\section{Appendix}

In this section, we will give the proofs for the two lemmas that have been used in proving Theorem~\ref{th3.2}. One is about the local characterization of linear maps between projective spaces and the other one is a characterization of linear isometric embeddings between two possibly degenerate indefinite Hermitian vector spaces.

\begin{lemma} \label{le3.1}
	Let \(U\) be a connected open subset of  $\mathbb{P}^n$. If $1\leq\ell'\leq \ell$ and \(h: U \rightarrow \mathbb{P}^m\) is a holomorphic map such that for  every $\ell$-dimensional linear subspace \(L \subset \mathbb{P}^n\), we have \(h(U\cap L) \subset L'\)  for some $\ell'$-dimensional linear subspace \(L'\subset\mathbb{P}^m\). Then, either the image of $h$ is contained in an $\ell'$-dimensional linear subspace or \(h\) extends  to a linear rational map \(\widehat{h}:\mathbb{P}^n \dasharrow \mathbb{P}^m \).
\end{lemma}
\begin{proof}
Without loss of generality, we may assume that $\ell'$ is the dimension of the linear span of $h(U\cap L)$ for a general $\ell$-dimensional linear subspace $L\subset\mathbb P^n$. 
Let $L_1,L_2\subset\mathbb P^n$ be a general choice of two $\ell$-dimensional linear subspaces  such that $\dim(L_1\cap L_2)=\ell-1$. Then, either $h(U\cap L_1)$ and $h(U\cap L_2)$ are contained in the same $\ell'$-dimensional linear subspace or $h(U\cap L_1\cap L_2)$ is contained in some $(\ell'-1)$-dimensional linear subspace. For the former possibility, since our choice is general, it follows that $h(U\cap L)$ is contained in the same $\ell'$-dimensional linear subspace for every $\ell$-dimensional linear subspace $L\subset\mathbb P^n$ and hence so is the entire image of $h$. Thus, by induction, we deduce that either the image of $h$ is contained in an $\ell'$-dimensional linear subspace or for every $(\ell-\ell'+1)$-dimensional linear subspace $\Phi\subset\mathbb P^n$, $h(U\cap\Phi)$ is contained in some line in $\mathbb P^m$. In the latter case, we have, in particular, for every line $\mathcal L\subset\mathbb P^n$, $h(U\cap\mathcal L)$ is contained in some line, which implies that  either the image of $h$ is contained in a line or $h$ extends to a linear rational map.  (See~\cite{ng2015holomorphic}, Lemma 4.1 therein, for a proof of this fact.)
\end{proof}

\begin{lemma}\label{linear algebra}
Let $V$, $W$ be two complex vector spaces equipped with indefinite (possibly degenerate) Hermitian form $\langle\cdot,\cdot\rangle_V$ and $\langle\cdot,\cdot\rangle_W$ respectively. If $G:V\rightarrow W$ is a linear map such that null vectors are mapped to null vectors, then there exists $\lambda\in\mathbb R$ such that $\langle Gu,Gv\rangle_W=\lambda\langle u,v\rangle_V$ for every $u,v\in V$.
\end{lemma}
\begin{proof}
Denote by $\|\cdot\|_V^2$ and $\|\cdot\|_W^2$ the corresponding norm squares for the Hermitian forms.
Choose an orthogonal decomposition $V=V_+\oplus V_-\oplus V_0$, where $V_+$, $V_-$ and $V_0$ are subspaces on which the restriction of $\langle\cdot,\cdot\rangle_V$ is positive definite, negative definite and zero respectively. Pick any $v_+\in V_+$, $v_-\in V_-$ and $v_0\in V_0$ such that $\|v_+\|_V^2=-\|v_-\|_V^2$, then for every $\alpha,\beta,\gamma\in\mathbb C$ such that $|\alpha|^2=|\beta|^2$, we have $\|\alpha v_++\beta v_-+\gamma v_0\|_V^2=0$ and thus $\|G(\alpha v_++\beta v_-+\gamma v_0)\|_W^2=0$ by hypotheses. By choosing $\gamma\in\mathbb C$ arbitrarily and $(\alpha,\beta)=(1,e^{i\theta})$, where $\theta\in\mathbb R$, we see that this implies  
$$
\langle Gv_+,Gv_-\rangle_W=\langle Gv_+,Gv_0\rangle_W=\langle Gv_-,Gv_0\rangle_W=0
$$
and
$$
\|Gv_+\|_W^2=-\|Gv_-\|_W^2.
$$

If we fix $v_-\in V_-$ and vary $v_+$, then from the second equation we deduce that there exists $\lambda\in\mathbb R$ such that $\|Ge_+\|_W^2=\lambda$ for every unit vector $e_+\in V_+$. Thus, $\|Gv_+\|_W^2=\lambda\|v_+\|_V^2$ for every $v_+\in V_+$. By polarization, we also get $\langle Gu_+,Gv_+\rangle_W=\lambda\langle u_+,v_+\rangle_V$, for every $u_+,v_+\in V_+$. Reversing the role of $V_+$ and $V_-$ and using the second equation above, we can deduce that with the same $\lambda$, we also have $\langle Gu_-,Gv_-\rangle_W=\lambda\langle u_-,v_-\rangle_V$, for every $u_-,v_-\in V_-$. In addition, $\langle Gu_0,Gv_0\rangle_W=\lambda\langle u_0,v_0\rangle_V$ holds trivially for every $u_0,v_0\in V_0$ by hypotheses. 

Finally, from the first equation, we see immediately that $G(V_+)\perp G(V_-)$, $G(V_+)\perp G(V_0)$ and $G(V_-)\perp G(V_0)$ and the desired result now follows.
\end{proof}

\bigskip
{\bf Acknowledgements.} The first author was partially supported by Science and Technology Commission of Shanghai Municipality (STCSM) (No. 13dz2260400).


\begin{thebibliography}{99}  

	\bibitem[BEH]{baouendi2011holomorphic} Baouendi, M. S., Ebenfelt, P., Huang, X.: Holomorphic mappings between hyperquadrics with small signature difference. {\it Amer. J. Math.} {\bf 133} (2011), 1633-1661.
		
	\bibitem[BH]{baouendi2005super} Baouendi, M. S., Huang, X.: Super-rigidity for holomorphic mappings between hyperquadrics with positive signature. {\it J. Diff. Geom.} {\bf 69} (2005), 379-398.
	
	\bibitem[EB]{ebenfelt} Ebenfelt, P.: Uniformly Levi-degenerate {CR}-manifolds: the $5$-dimensional case. {\it Duke Math. J.} {\bf 110} (2001), 37-80.	
	
	
	\bibitem[GN]{gao2018rational} Gao, Y., Ng, S.-C.: On rational proper mappings among generalized complex balls. {\it Asian J. Math.} {\bf 22} (2018), 355-380.
	
	\bibitem[GM]{gm} Gaussier, H.,  Merker, J.: A new example of a uniformly Levi degenerate hypersurface in {$\mathbb C^3$}. {\it Ark. Mat.} {\bf 41} (2003), 85-94.

	\bibitem[IK]{ik} Isaev, A. V., Kossovskiy, I. G.: Continuation of CR-automorphisms of Levi degenerate hyperquadrics to the projective space. {\it Illinois J. Math.} {\bf 54} (2010), 747-752.
	
\bibitem[KLX]{klx} Kossovskiy, I. G., Lamel, B., Xiao, M.: Regularity of CR-mappings into Levi-degenerate hypersurfaces. {\it to appear in Comm. Anal. Geom.} 

	
    \bibitem[N1]{ng2010holomorphic}  Ng, S.-C.: On holomorphic isometric embeddings of the unit disk into polydisks. {\it Proc. Amer. Math. Soc.} {\bf 138} (2010), 2907-2922.

    \bibitem[N2]{Ng1}  Ng, S.-C.: Cycle spaces of flag domains on Grassmannians and rigidity of holomorphic mappings. {\it Math. Res. Lett.} {\bf 19} (2012), 1219-1236.

	 \bibitem[N3]{ng2013proper}  Ng, S.-C.: Proper holomorphic mappings on flag domains of {$ SU(p, q) $} type on projective spaces. {\it Michigan Math. J.} {\bf 62} (2013), 769-777.
	 
	 \bibitem[N4]{ng2015holomorphic}  Ng, S.-C.: Holomorphic double fibration and the mapping problems of classical domains. {\it Int. Math. Res. Not.} {\bf 2015} (2015), 291-324.
	 
	 
	 \bibitem[SE]{seo}  Seo, A.: New examples of proper holomorphic maps among symmetric domains. {\it Michigan Math. J.} {\bf 64} (2015), 435-448.
	
	 

	\end{thebibliography}
\end{document}